\documentclass[12pt]{elsarticle}

\usepackage{amsmath,amssymb,amsthm,amscd}
\usepackage{multicol,multirow}

\makeatletter
\newtheorem*{rep@theorem}{\rep@title}
\newcommand{\newreptheorem}[2]{%
\newenvironment{rep#1}[1]{%
 \def\rep@title{#2 \ref{##1}}%
 \begin{rep@theorem}}%
 {\end{rep@theorem}}}
\makeatother

\newtheorem{lemma}{Lemma}[section]

\newtheorem{theorem}[lemma]{Theorem}

\newreptheorem{thma}{Theorem}

\theoremstyle{definition}
\newtheorem{definition}[lemma]{Definition}
\newtheorem{example}[lemma]{Example}
\newtheorem{remark}[lemma]{Remark}
\newtheorem{conjecture}[lemma]{Conjecture}
\newtheorem{corollary}[lemma]{Corollary}
\newtheorem*{cor-nonum}{Corollary}

\numberwithin{equation}{section}

\def\det{\text{det} \thinspace}

\newcommand{\Z}{{\mathbb Z}}
\newcommand{\C}{{\mathbb C}}
\newcommand{\R}{{\mathbb R}}

\newcommand{\Q}{{\mathbb Q}}
\newcommand{\PP}{{\mathbb P}}
\newcommand{\RP}{{\mathbb {RP}}}
\newcommand{\CP}{{\mathbb {CP}}}
\newcommand{\proofend}{\hfill$\Box$\bigskip}

\def\proof{\paragraph{Proof}}

\journal{Journal of Geometry and Physics}
\begin{document}

\begin{frontmatter}



\title{Non-integrability vs. integrability in pentagram maps}


\author{Boris Khesin}

\address{Department of Mathematics, University of Toronto, Toronto, ON M5S 2E4, Canada }

\author{Fedor Soloviev}

\address{The Fields Institute, Toronto, ON and CRM, Montreal, QC, Canada }

\begin{abstract}
We revisit recent results on integrable cases for higher-dimensional generalizations of the 2D
pentagram map: short-diagonal, dented, deep-dented, and corrugated versions, and define a universal class of pentagram maps, which are proved to possess projective duality.
We show that in many cases the pentagram map cannot be included into  integrable flows
as a time-one map, and discuss how the corresponding notion of discrete integrability
can be extended to include jumps between invariant tori.
We also present a numerical evidence that certain
generalizations of the integrable 2D pentagram map are non-integrable and  present a conjecture for  a necessary condition of their discrete integrability.
\end{abstract}


\end{frontmatter}



~

\bigskip
~

\bigskip

The goal of this paper is three-fold. First we revisit the recent progress in finding integrable
generalizations of the 2D pentagram map. Secondly, we discuss a natural framework
for a notion of a discrete integrable Hamiltonian map.
It turns out that the Arnold--Liouville theorem on existence of invariant tori admits a
natural generalization to allow discrete dynamics  with jumps between invariant tori,
which is relevant for many pentagram maps.
Lastly, we define a universal class of pentagram-type maps, describe a projective duality for them, and present numerical evidence for non-integrability of
several pentagram maps in 2D and 3D. In view of many new integrable generalizations found recently, a search for a non-integrable  generalization of the pentagram map
was brought into light, and the
examples presented below might help focusing the efforts for  such a search.

\section{Types of pentagram maps}\label{intro}

Recall that the pentagram map is a map on plane convex
polygons considered up to their projective equivalence,
where a new polygon is spanned by the shortest diagonals of the initial one, see  \cite{Schwartz}.
It exhibits quasi-periodic behaviour of (projective classes of) polygons in 2D under iterations, which indicates hidden integrability. The integrability of this map was proved in  \cite{OST99}, see also \cite{FS}.

While the pentagram map is in a sense unique in 2D, its generalizations
to higher dimensions allow more freedom.
It turns out that while there seems to be  no natural generalization of this map to polyhedra,
 one can suggest several natural integrable extensions of the pentagram map
to the space of generic twisted polygons in higher dimensions.

\smallskip

\begin{definition}
A {\it twisted $n$-gon} in a projective space $\PP^d$ with a monodromy $M \in SL_{d+1}$
is a doubly-infinite sequence of points $v_k\in \PP^d$, $k\in \mathbb Z$,
such that
$v_{k+n} =  M \circ v_k$ for each $k\in \mathbb Z$, and where $M$ acts naturally on $\PP^d$.
We assume that the vertices $v_k$ are in general position (i.e., no $d+1$ consecutive vertices lie in the same hyperplane in $\PP^d$),
and denote by   ${\mathcal P}_n$ the space of generic twisted $n$-gons considered up to the projective equivalence.
\end{definition}

We use projective spaces defined over reals $\R$ (as the easiest ones to visualize), over complex numbers $\C$ (to describe algebraic-geometric
integrability), and over rational numbers $\Q$ (to perform a non-integrability test).
All definitions below work for any base field.
General pentagram  maps are defined as follows.

\smallskip

\begin{definition}\label{def:I-diag}
We define 3 types of  diagonal hyperplanes for a given twisted polygon $(v_k)$ in $\PP^d$.
$a)$ The {\it short-diagonal hyperplane} $P^{\rm sh}_k$ is defined as the  hyperplane passing through $d$ vertices of the $n$-gon by taking every other vertex starting with $v_k$:
$$
P^{\rm sh}_k:=(v_k, v_{k+2}, v_{k+4},..., v_{k+2(d-1)})\,.
$$
$b)$ The {\it dented  diagonal plane hyperplane} $P^m_k$ for a fixed $m=1,2,...,d-1$
is the hyperplane passing through all vertices
from $v_k$ to $v_{k+d}$ but one, by skipping only the vertex $v_{k+m}$:
$$
P^m_k:=(v_k, v_{k+1},...,v_{k+m-1},v_{k+m+1},v_{k+m+2},...,v_{k+d})\,.
$$
$c)$ The {\it deep-dented diagonal plane hyperplane} $P^m_k$ for fixed positive integers $m$ and $p\ge 2$ 
is the hyperplane as above that passes through consecutive vertices, except for one jump, when it skips $p-1$ vertices $v_{k+m}$,...,$v_{k+m+p-2}$:
$$
P^{m,p}_k:=(v_k, v_{k+1},...,v_{k+m-1},v_{k+m+p-1},v_{k+m+p},...,v_{k+d+p-2})\,.
$$
(Here $P_k^{m,2}$ corresponds to $P_k^m$ in $b)$.)
Now the corresponding {\it pentagram maps} $T_{\rm sh}$, $T_m$, and $T_{m,p}$
are defined on generic twisted polygons  $(v_k)$ in $\PP^d$ by intersecting $d$ consecutive diagonal hyperplanes:
$$
T v_k:=P_{k}\cap P_{k+1}\cap ...\cap P_{k+d-1}\,,
$$
where each of the maps  $T_{\rm sh}$, $T_m$, and $T_{m,p}$ uses the
definition of the corresponding hyperplanes  $P^{\rm sh}_k$, $P^m_k$, and $P^{m,p}_k$.
 These pentagram maps are generically defined on the classes of projective equivalence of twisted polygons $T:{\mathcal P}_n\to {\mathcal P}_n$.
\end{definition}

\begin{example}
For $d=2$ one can have only $m=1$ and the definitions of $T_{\rm sh}$ and $T_m$
coincide  with the standard 2D pentagram map $T_{\rm st}$ in \cite{Schwartz}
(up to a shift in vertex numbering).
The deep-dented maps  $T_{1,p}$ in 2D are the maps
$T_{1,p} v_k:=(v_k,v_{k+p}) \cap (v_{k+1},v_{k+p+1})$ obtained by intersecting deeper diagonals
of  twisted polygons, see Figure \ref{fig:T13_2D}.
\begin{figure}[hbtp!]
\centering
\includegraphics[width=2.8in]{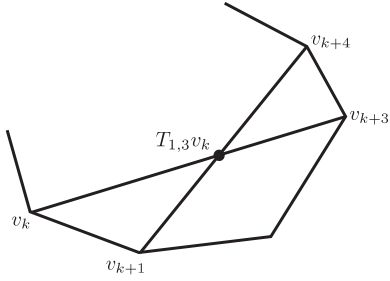}
\caption{\small Deeper pentagram map $T_{1,3}$ in 2D.} \label{fig:T13_2D}
\end{figure}

For $d=3$ the map $T_{\rm sh}$ uses
the diagonal planes $P^{\rm sh}_k:=(v_k, v_{k+2}, v_{k+4})$,
while for the dented maps $T_1$ and $T_2$ one has
$P^1_k=(v_k, v_{k+2}, v_{k+3})$ and $P^2_k=(v_k, v_{k+1}, v_{k+3})$ respectively, see Figure \ref{fig:T-3D}.
\begin{figure}[hbtp!]
\centering
\includegraphics[width=5.2in]{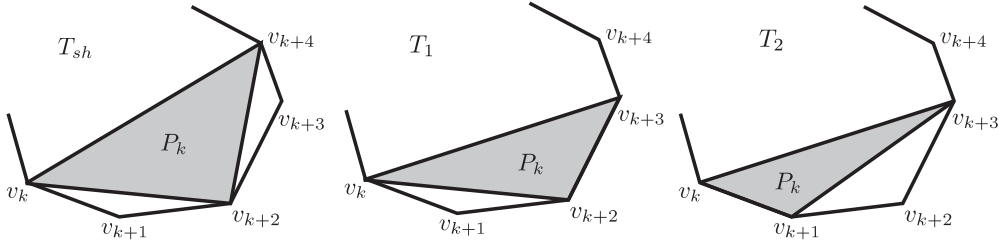}
\caption{\small Different diagonal planes in 3D: for $T_{\rm sh}, T_1,$ and $ T_2$.}
\label{fig:T-3D}
\end{figure}
\end{example}

\begin{theorem}\label{th-lax}
The short-diagonal $T_{\rm sh}$, dented  $T_m$  and deep-dented $T_{m,p}$ maps are integrable in any dimension $d$ on both twisted and closed
$n$-gons in a sense that they admit Lax representations with a spectral parameter.
\end{theorem}

The integrability
of the standard 2D pentagram map $T_{\rm st}:=T_{\rm sh}=T_m$ was proved in  \cite{OST99},
while its Lax representation was found in \cite{FS}.
In \cite{GSTV} integrability of the pentagram map for corrugated polygons (which we discuss below) was proved, which implies integrability of the  maps $T_{1,p}$ in 2D.

For higher pentagram maps in any dimension $d$, their Lax representations
with a spectral parameter were found in \cite{KS}. The dependence on spectral parameter
was based on the scale invariance of such maps,
which was proved in \cite{KS} for 3D, and  in \cite{Beffa_scale} for higher $d$.
For the dented and deep-dented pentagram maps their  Lax representations
and scale invariance in any dimension $d$ were established in \cite{KS13}.
We present  formulas for those Lax representations in Section \ref{laxrep}.
Such representations with a spectral parameter provide first integrals of the maps
(as the coefficients
of the corresponding spectral curves) and allow one to use algebraic-geometric
machinery to prove various integrability properties.

In \cite{KS, KS13} we proved that the proposed Lax representation implies algebraic-geometric
integrability for the maps $T_{\rm sh}$, $T_1$, $T_2$ in 3D. In particular, this means that
the space of twisted $n$-gons in the complex space ${\mathbb{CP}}^3$
is generically fibered into (Zariski open subsets of) tori whose dimension is described in terms of
$n$.
In Section \ref{general} we discuss  features of the pentagram  maps which emphasize their discrete nature.

\smallskip

\begin{definition}\label{def:IJ}
More generally, one can define {\it generalized pentagram maps} $T_{I,J}$ on (projective equivalence classes of)
twisted polygons in $\PP^d$, associated with a $(d-1)$-tuple of integers $I$ and $J$:
the {\it jump tuple} $I=(i_1, ..., i_{d-1})$ determines which vertices define the diagonal hyperplanes $P^I_k$:
\begin{equation}\label{eq:P^I}
P_k^I:=(v_k, v_{k+i_1},...,v_{k+i_1+...+i_{d-1}}),
\end{equation}
while the {\it intersection tuple} $J=(j_1, ..., j_{d-1})$ determines which of the hyperplanes
to intersect in order to get the image of the  point $v_k$:
$$
T_{I,J}v_k:=P^I_k\cap P^I_{k+j_1}\cap...\cap P^I_{k+j_1+...+j_{d-1}}\,.
$$
\end{definition}

In general, the integrability of $T_{I,J}$ is yet unknown, but there exists the following duality between such pentagram maps:
\begin{equation}\label{eq:duality}
T_{I,J}^{-1}=T_{J^*,I^*}\circ Sh\,,
\end{equation}
where $I^*$ and $J^*$ stand for the   $(d-1)$-tuples  taken in the opposite order and  $Sh$
is any shift in the indices  of polygon vertices, see \cite{KS13}.
In particular, the maps $T_{I,J}$ and $T_{J^*,I^*}$ are integrable or non-integrable simultaneously.

The pentagram maps $T_{\rm sh}$,  $T_m$,  and $T_{m,p}$ considered above
correspond to  $J=(1,...,1)$ (cf. Definitions \ref{def:I-diag} and
\ref{def:IJ}). The duality  \eqref{eq:duality} of $T_{I,J}$ and
$T_{J^*,I^*}$ along with Theorem \ref{th-lax} imply  integrability
of the maps with $I=(1,...,1)$ and appropriate $J$'s.

\medskip
The simplest pentagram map which is neither short-diagonal, nor dented or deep-dented
appears in dimension $d=3$ and corresponds to $I=(2,3)$ and $J=(1,1)$.
We conjecture that it is indeed non-integrable and outline supporting evidence from computer
experiments in Sections \ref{nonint2D} and \ref{nonint3D}, along with several other cases both integrable and not.

\begin{remark} In \cite{KS, KS13} it was also proved that the continuous limit
of any short-diagonal or
dented pentagram map (and more generally, of any generalized pentagram map) in
${\mathbb{RP}}^d$  is the $(2, d+1)$-KdV flow of the Adler-Gelfand-Dickey  hierarchy on the circle.
For 2D this is the classical  Boussinesq equation  on the circle: $u_{tt}+2(u^2)_{xx}+u_{xxxx}=0$,
which appears as the continuous limit of the 2D pentagram map \cite{OST99}.
\end{remark}

\begin{remark} Note also  that a  different integrable generalization to higher dimensions was proposed in \cite{GSTV}, where the pentagram map was defined not on generic, but on the so-called
{\it  corrugated} polygons. These are  twisted polygons in
$\PP^d$, whose vertices $v_{k-1}, v_{k}, v_{k+d-1},$ and $v_{k+d}$
span a projective two-dimensional plane  for every $k\in \mathbb Z$. The pentagram map
$T_{\rm cor}v_k:=(v_{k-1}, v_{k+d-1})\cap(v_{k}, v_{k+d})$
on corrugated polygons  turns out to be integrable  and  admits an explicit description of the Poisson structure, a cluster algebra structure, and other  interesting features \cite{GSTV}.
Furthermore, it turns out that the pentagram map $T_{\rm cor}$
can be viewed as a particular case of the dented pentagram map:
\end{remark}


\begin{theorem} {\bf (\cite{KS13})}
{\it This pentagram map $T_{\rm cor}$ is a restriction of the dented pentagram map $T_m$ for any $m=1,..., d-1$ from generic n-gons ${\mathcal P}_n$ in $\PP^d$ to corrugated ones
 (or differs from it by a shift in vertex indices). In particular, these restrictions for different $m$  coincide modulo an index shift.}
\end{theorem}


\section{Formulas for Lax representations}\label{laxrep}

In this section we recall explicit formulas of the Lax representation for pentagram maps.
First we introduce coordinates on the space ${\mathcal P}_n$ of generic
twisted $n$-gons in $\PP^d$  considered over $\C$.
For simplicity, we focus only on the case when $gcd(n,d+1)=1$ (see the general case
in \cite{KS, KS13}).

One can show that there exists a lift of the vertices $v_k=\phi(k) \in \CP^d$
to the vectors $V_k \in \C^{d+1}$ satisfying $\det(V_j, V_{j+1}, ..., V_{j+d})=1$
and $ V_{j+n}=MV_j,\; j \in \Z,$ where $M\in SL_{d+1}(\C)$.
(Strictly speaking, this lift is not unique, because it is defined up to a simultaneous
multiplication of all vectors by $(-1)^{1/(d+1)}$,
but  the coordinates introduced below have the same values for all lifts.)\footnote{Note also
that over $\R$  for odd $d$  to obtain the lifts of $n$-gons from $\RP^d$ to $\R^{d+1}$
one might need to switch the sign of the monodromy matrix: $M \to -M \in SL_{d+1}(\R)$, since the field is not algebraically closed. These monodromies correspond to the same projective monodromy
in $PSL_{d+1}(\R)$.}

The coefficients of the following difference equation
$$
V_{j+d+1} = a_{j,d} V_{j+d} + a_{j,d-1} V_{j+d-1} +...+ a_{j,1} V_{j+1} +(-1)^{d} V_j,\quad j \in \Z.
$$
turn out to be $n$-periodic in $j$ due to the monodromy relation on vectors $V_j$, and, in particular, coefficients $\{a_{j,k}~|~j=0,...,n-1, k=1,...,d\}$
play the role of the coordinates on the space ${\mathcal P}_n$.
The dimension of the space $ {\mathcal P}_n$ of generic $n$-gons in $\CP^d$ is
$\text{dim} \thinspace {\mathcal P}_n=nd$.

Now we are in a position to define Lax representations for the maps $T_{\rm sh}$ and $T_m$.
The above pentagram maps can be presented in the Lax form
$L_{j,t+1}(\lambda) = P_{j+1,t}(\lambda) L_{j,t}(\lambda) P_{j,t}^{-1}(\lambda)$
for an appropriate matrix $P_{j,t}(\lambda)$, where $\lambda$ is a spectral parameter.
We present the $L$ matrices below, while the explicit expression for $P_{j,t}(\lambda)$ is complicated and is not required for our analysis.\footnote{One can recover the $P$-matrix from the coordinate formulas of the map and the fact that the ordered product $L_{n-1}...L_1 L_0$ of $L$-matrices transforms by conjugation
after the application of the map, see  \cite{KS} for more detail.} The pentagram map corresponds to the time evolution
$t\to t+1$ in the Lax matrix.

\begin{theorem} {\bf (=Theorem \ref{th-lax}$'$)}
Lax representations with a spectral parameter for the above pentagram maps are given by the following $L$-matrices:
\[
L_{j,t}(\lambda) =
\left(
\begin{array}{cccc|c}
0 & 0 & \cdots & 0    &(-1)^d\\ \cline{1-5}
\multicolumn{4}{c|}{\multirow{4}*{$D(\lambda)$}} & a_{j,1}\\
&&&& a_{j,2}\\
&&&& \cdots\\
&&&& a_{j,d}\\
\end{array}
\right)^{-1},
\] 
where $D(\lambda)$ is the following diagonal $(d \times d)$-matrix:
\begin{itemize}
\item for the map $T_{\rm sh}$, $D(\lambda)={\rm diag}(\lambda,1,\lambda,1,...,1,\lambda)$ for odd $d$ and $D(\lambda)={\rm diag}(1,\lambda,1,...,1,\lambda)$ for even $d$;
\item for the map $T_m$, $D(\lambda)={\rm diag}(1,...,1,\lambda, 1,...1)$, where
the spectral parameter $\lambda$ is situated at the $(m+1)${\rm th} place.
\end{itemize}

\end{theorem}

The construction of a Lax representation for the deep-dented  pentagram map $T_{m,p}$
relies on lifting generic polygons from $\CP^d$ to so-called partially  corrugated polygons
in a bigger space $\CP^{d+p-2}$. Then the corresponding Lax representation
for the deep-dented maps can be obtained from the above Lax form for map on polygons
in $\CP^{d+p-2}$ by restricting it to the subset of partially  corrugated ones,
see details in \cite{KS13}.


\section{Discrete integrability}\label{general}

The above Lax representation allows one to give a more detailed description of the  dynamics. In particular, in lower dimensions one can explicitly express the pentagram maps as a discrete dynamics on the Jacobian of the corresponding spectral curve. The following result is a corollary of that description.

\begin{theorem}
The above integrable pentagram maps on twisted $n$-gons in $\CP^d$ cannot be included into a
 Hamiltonian  flow as its time-one map, at least for some values of $n, m,$ and $d$.
\end{theorem}

\proof
In \cite{FS, KS, KS13} we gave a detailed description of the pentagram maps $T_{\rm sh}$ and
$T_m$ in 2D and 3D (denote these maps by $T_*$).
It turned out that for even $n$ one observes the staircase-like dynamics
on the Jacobian of the corresponding spectral curve. In the space $\mathcal P_n$
of generic twisted $n$-gons this corresponds to the following phenomenon.
This space is a.e. fibered into (Zariski open subsets of) complex tori, which are invariant for
 the square  $T_*^2$ of  the pentagram map, but not for the map $T_*$ itself.
(More generally, the tori are invariant for a certain power $T_*^q$,
while we set $q=2$ for the rest of the proof.)
 In turn, the map $T_*$ sends almost every $n$-gon from the space $\mathcal P_n$
 to jump between two tori. The square of this map is a shift on each torus.

Now assume that such a map $T_*$ were the time-one map of a smooth autonomous
Hamiltonian field $v$ on  $\mathcal P_n$.  Then this Hamiltonian field
admits the same fibration a.e. into invariant tori, since $T_*^2$ is its time-two map and its frequencies are known to be nondegenerate. Then the flow of this field $v$
would describe the linear evolution on tori, and hence it would be integrable itself. The map
$T_*$ is by assumption the time-one map of the same flow, and hence it must have the same
invariant tori , rather than
jumping between them. This contradiction proves that inclusion into a flow is impossible.
\proofend

Note that  the dynamics of (partially) corrugated polygons, described in \cite{KS13},
allows jumps between 3 different  tori for some values of $n, m, p$, and $d$.

We conjecture that the pentagram dynamics cannot be included into a flow for {\it all values} of
$n, m,$ and $d$ (even when the above simple argument does not already work).
The consideration and examples above suggest the following generalization of a discrete
integrable Hamiltonian system. It can be regarded as a particular case
of an integrable correspondence \cite{Ves}.

\begin{definition}
Suppose that $(M,\omega)$ is a $2n$-dimensional symplectic manifold and $I_1,...,I_n$ are $n$ independent functions in involution.
Let $M_{\bf c}$ be a  (possibly disconnected) level set of these functions:
$
M_{\bf c}=\{ x\in M~|~ I_j(x)=c_j,\; 1 \le j \le n \}.
$
A map $T: M \to M$ is called {\it generalized integrable} if
\begin{itemize}
\item it is symplectic, i.e., $T^*\omega=\omega$;
\item it preserves the integrals of motion: $T^*I_j  \equiv I_j,\; 1 \le j \le n$;
\item there exists a positive  integer $q\ge 1$ such that the map $T^q$ leaves all connected components of level sets $M_{\bf c}$ invariant for all ${\bf c}=(c_1,...,c_n)$.
\end{itemize}
\end{definition}

In other words, the  $q$th iteration $T^q$ of the map $T$ is integrable in the usual sense.
This definition lists almost verbatim the assertions of the Arnold--Liouville theorem \cite{A89}
for continuous flows, which implies that one has a conditionally periodic motion for the map $T^q$ and its ``integrability by quadratures.''
The difference with the classical case, corresponding to $q=1$, may occur if  level sets $M_{\bf c}$
are disconnected, since the discrete map can ``jump'' from one component to another.

Our analysis (in the complexified case) shows that the pentagram maps
$T_{\rm sh}$ for $d=3$ and even $n$ (in that case $q=2$) and $T_{\rm cor}$ in the corrugated case for $d=3$ and $n=6l+3$ (then one has $q=3$) are generalized integrable.
Note that compact connected components of generic level sets  $M_{\bf c}$ are tori, and
the map $T$ can be used to establish an isomorphism of different connected components.
 Under such an isomorphism one obtains a ``staircase'' dynamics on the same torus
 (as discussed, e.g., in Theorems B in \cite{FS, KS}).


\section{Universal pentagram maps}\label{sect:univ}

In this section we define a general class of pentagram maps
in any dimension, which allows one to intersect different diagonals at
each step.

\begin{definition}\label{def:sub-univ}
Let $(v_k)$ be a generic twisted $n$-gon in $\PP^d$. We fix $d$, pairwise
different, $d$-tuples $I_1, ..., I_d$ of integers which are jump tuples defining  $d$ hyperplanes $P_k^{I_1}, ..., P_k^{I_d}$, i.e.,  each
hyperplane $P_k^{I_\ell}$ passes through the vertices defined by its own
jump $d$-tuple $I_\ell=(i_{\ell ,1},..., i_{\ell ,d}) $:
$$
P^{I_\ell}_k:=(v_{k+i_{\ell ,1}}, v_{k+i_{\ell ,2}},...,v_{k+i_{\ell ,d}})\,.
$$
Now we define the {\it skew pentagram map} $T_{I_1,...,I_d}:\mathcal P_n \to \mathcal P_n$, where
vertices of a new $n$-gon are obtained  by intersecting these $d$ hyperplanes  $P_k^{I_\ell}, \ell=1,...,d$:
$$
T_{I_1,...,I_d}v_k:=P^{I_1}_k\cap P^{I_2}_{k}\cap...\cap P^{I_d}_{k}.
$$
\end{definition}

\begin{remark}
$a)$ The general pentagram map $T_{I,J}$ described in Definition  \ref{def:IJ}  is a particular
case of the skew pentagram map $T_{I_1,...,I_d}$: one can obtain both the jump $(d-1)$-tuple $I$
and the intersection $(d-1)$-tuple $J$ from the set of $d$-tuples $I_1, ..., I_d$.

$b)$ The pentagram maps in \cite{Beffa}  are defined by intersecting
a segment $(v_{k-1}, v_{k+1})$ with $P_k^I$ for an appropriate choice of jumps $I$. This
is a particular case of $T_{I_1,...,I_d}$   with $I_d=I$, and $d$-tuples $I_1,...,I_{d-1}$ all containing $(-1, +1,  ...)$, so that the planes
$P_k^{I_\ell}$  all contained the pair of vertices $(v_{k-1}, v_{k+1})$, while their other
vertices were different. In this case
$P^{I_1}_k\cap P^{I_2}_{k}\cap...\cap P^{I_{d-1}}_{k}=(v_{k-1},v_{k+1}).$

In particular, the class of pentagram maps $T_{I_1,...,I_d}$ contains
pentagram maps defined by taking intersections of
subspaces of complimentary dimensions (and spanned by vertices
$(v_k)$) to obtain a point as an intersection. For instance, the map defined by
the intersection $(v_k,v_{k+3})\cap(v_{k+1},v_{k+2},v_{k+4})$ of a segment and a plane can be equivalently defined as the intersection of three planes:
$(v_k,v_{k+3}, v_{k+5}) \cap(v_k,v_{k+3}, v_{k+6})\cap(v_{k+1},v_{k+2},v_{k+4})$.
Note that the intersections of
hyperplanes provide a more general definition, since their intersection subspaces might not necessarily be spanned by vertices $(v_k)$ themselves, but by their linear combinations.
 \end{remark}

Finally, define a universal pentagram map by starting with $d$
polygons.

\begin{definition}\label{def:univ}
Let $(v_k^\ell)$ be $d$ twisted polygons in $\PP^d$, $\ell=1,...,d$
and $k\in \Z$, with the same monodromy matrix $M\in SL_{d+1}$.
Now we fix two sets of $d$-tuples, jump tuples $I_1, ..., I_d$
and intersection tuples $J_1,..., J_d$. Let $I_\ell=(i_{\ell ,1},..., i_{\ell ,d}) $ and
$J_p=(j_{p,1},..., j_{p,d})$.
Define the hyperplane
$$
P^{I_\ell}_k:=(v_{k+i_{\ell ,1}}^1, v_{k+i_{\ell ,2}}^2,...,v^d_{k+i_{\ell ,d}})\,.
$$
i.e., this plane $P^{I_\ell}_k$ uses one vertex from each $n$-gon. Now one can define
$d$ skew pentagram maps, or rather a {\it universal pentagram map}, whose image
consists of $d$ twisted $n$-gons: for every $p=1,...,d$ the map $T_p$ uses the
corresponding intersection tuple $J_p$:
$$
T_pv_k:=P^{I_1}_{k+j_{p ,1}}\cap P^{I_2}_{k+j_{p,2}}\cap...\cap P^{I_d}_{k+j_{p ,d}}.
$$
Thus one obtains a universal map $T_{\mathcal I,\mathcal J}$ on $d$-tuples of twisted $n$-gons
which is associated with two sets of tuples $\mathcal I=(I_1,...,I_d)^t$ and
$\mathcal J=(J_1,..., J_d)^t$, where these sets $\mathcal I=(i_{\ell ,s})$ and
$\mathcal J=(j_{p,s})$ can be thought of as two $(d\times d)$-matrices composed of $d$-tuples
$I_1,...,I_d$ and, respectively, $J_1,...,J_d$ written as their rows.
\end{definition}

\begin{theorem} \label{thm:duality}
Universal pentagram maps possess the following duality:
$$
T^{-1}_{\mathcal I,\mathcal J}=T_{-\mathcal I^*,-\mathcal J^*},
$$
where $\mathcal I^*$ and $\mathcal J^*$ stand for the transposed matrices $\mathcal I^*=(i_{s,\ell })$ and $\mathcal J^*=(j_{s,p})$ respectively.
\end{theorem}

This duality generalizes the one \eqref{eq:duality} for the $T_{I,J}$ pentagram maps.
Note that  the description in terms of jump $(d-1)$-tuples does not acquire minuses, but gives a similar equality only up to a shift of indices, while in the description in terms of  $d$-tuples
one fully specifies the indices.

\proof
To prove this theorem we modify  the notion of a duality map, cf. \cite{OST99, KS13}.

\begin{definition}
{\rm
Given $d$ generic sequences of points $\phi_\ell(j) \in \RP^d, \; j \in \Z,\; \ell=1,...,d$  and
a $d$-tuple $I=(i_1,...,i_{d})$
we define the following {\it sequence of hyperplanes} in $\RP^d$ enumerated by $j$:
$$
\alpha_I(\phi_\star(j)):=(\phi_1(j+i_1), \phi_2(j+i_2),..., \phi_d(j+i_{d}))\,,
$$
which is regarded as a sequence of points in the dual space: $\alpha_I(\phi_\star(j))\in (\RP^d)^*$,
$ j \in \Z$.
For $d$ tuples $I_1, ..., I_d$ we get $d$ sequences
$\alpha_{\mathcal I}(\phi_\star(j))=(\alpha_{I_1}(\phi_\star(j)),...,\alpha_{I_d}(\phi_\star(j)))
\in (\RP^d)^*\times ... \times  (\RP^d)^*$ enumerated by $j$. In other words, starting with $d$
sequences of points in $\RP^d$,  the map $\alpha_{\mathcal I}$ gives $d$
sequences of points in the dual space $(\RP^d)^*$.
}
\end{definition}

The universal pentagram map $T_{\mathcal I,\mathcal J}$ on $d$ twisted polygons can be defined as a composition of two such maps:
$T_{\mathcal I,\mathcal J}=\alpha_{\mathcal I}\circ\alpha_{\mathcal J}$.
Note that these maps by definition possess the following duality property:
$T_{\mathcal I,-\mathcal I^*}=\alpha_{\mathcal I}\circ \alpha_{-{\mathcal I}^*}=Id$.

For instance, in $\RP^2$ consider two twisted polygons $(v^1_k)$ and $(v^2_k)$ and
two 2-tuples $I=(1, -2)$ and $I_2=(-5,3)$.
Then in the dual space, by applying $\alpha_{\mathcal I}$ we obtain two twisted polygons, formed by lines
$P^1_k:=(v^1_{k+1}, v^2_{k-2})$ and $P^2_k:=(v^1_{k-5}, v^2_{k+3})$.
Then the vertex $v^1_k$ can be recovered from
$v^1_k=P^1_{k-1}\cap    P^2_{k+5}$, while $v^2_k=P^1_{k+2}\cap    P^2_{k-3}$, i.e. by applying the map  $\alpha_{-\mathcal I^*}$ to the sequences  $(P^1_k)$ and $(P^2_k)$. Similarly this works in any dimension.

Now we see that
$$
T_{\mathcal I,\mathcal J}\circ T_{-\mathcal J^*,-\mathcal I^*}
=\alpha_{\mathcal I}\circ \alpha_{\mathcal J}\circ \alpha_{-\mathcal J^*}\circ \alpha_{-{\mathcal I}^*}=Id\,,
$$
as required.
\proofend

Note that if $\mathcal J=-\mathcal J^*$, i.e. the matrix  $\mathcal J$ is skew-symmetric,
then the map $\alpha_{\mathcal J}$ is an involution: $\alpha_{\mathcal J}\circ \alpha_{\mathcal J}=Id$.

\begin{corollary}\label{cor:dual}
If $\mathcal J$ is skew-symmetric, then  the pentagram maps $ T_{\mathcal I,\mathcal J}$ and $T_{\mathcal J,\mathcal I}$   are conjugated  to each other, i.e.,
the map  $\alpha_{\mathcal J}$ takes the map $T_{\mathcal I,\mathcal J}$ on $d$-tuples of twisted $n$-gons in $\RP^d$ into the map
 $T_{\mathcal J,\mathcal I}$ on  $d$-tuples of twisted $n$-gons in $(\RP^d)^*$.
 In particular, all four maps $ T_{\mathcal I,\mathcal J}, T_{-\mathcal I^*,\mathcal J},
 T_{\mathcal J,\mathcal I}$ and $T_{\mathcal J,-\mathcal I^*}$ are integrable or non-integrable simultaneously.
\end{corollary}

\proof
First note that
$$
\alpha_{\mathcal J} \circ T_{\mathcal I,\mathcal J}\circ  \alpha^{-1}_{\mathcal J}
=\alpha_{\mathcal J} \circ (\alpha_{\mathcal I}\circ \alpha_{\mathcal J})\circ \alpha_{\mathcal J}
=\alpha_{\mathcal J} \circ \alpha_{\mathcal I}=T_{\mathcal J,\mathcal I}\,.
$$
Hence the pentagram map $T_{\mathcal I,\mathcal J}$ is conjugated to $T_{\mathcal J,\mathcal I}$. Furthermore, the pentagram maps
$T_{\mathcal I,\mathcal J}$ and $T_{\mathcal J,-\mathcal I^*}$, as well as
$ T_{\mathcal J,\mathcal I}$ and $T_{-\mathcal I^*,\mathcal J}$,
 are inverses to each other
for $\mathcal J=-\mathcal J^*$,
as  follows from Theorem \ref{thm:duality}.
This proves the corollary.
\proofend

\begin{conjecture}\label{conj:uni}
$a)$ All  universal  pentagram maps $T_{\mathcal I,\mathcal J}$ are discrete Hamiltonian systems
(i.e., preserve a certain Poisson structure), although not necessarily integrable.

$b)$  A necessary condition for integrability of the universal
pentagram maps $T_{\mathcal I,\mathcal J}$ is their equivalence to a map $T_{I,J}$,
see Definition \ref{def:IJ}
\end{conjecture}

In the next two sections we provide a numerical evidence to the  Conjecture~\ref{conj:uni}~$b)$ and explain why this equivalence to an appropriate  map $T_{I,J}$ cannot be sufficient for integrability.


\section{Non-integrability in 2D}\label{nonint2D}

{\it The classical case.}
In this section we are going to compare several pentagram maps in 2D.
To detect integrability we use the height criterion following \cite{H05} (see more references on `height' in \cite{GHRV09}).
Recall that the {\it height} of a rational number $a/b \in \Q$, written in the lowest terms, is $ht(a/b)=\text{max}(|a|,|b|)$.
We employ the cross-ratio coordinates $(x,y)$ (defined in \cite{OST99}) on the space of twisted $n$-gons ${\mathcal P}_n$ sitting inside $\Q\PP^2$ (i.e., having only rational values of coordinates).

\begin{definition}
{\rm
The {\it height} of a twisted $n$-gon $P \in {\mathcal P}_n$ in
$\Q\PP^2$ is defined as
$$
H(P):=\underset{0 \le i \le n-1}{\text{max}} \; \text{max}(ht(x_i),ht(y_i)).
$$
}
\end{definition}

We trace how fast the height of an initial $n$-gon grows with the
number of iterates of the pentagram map
(i.e., with an integer parameter $t$). We perform the comparison for
$n$-gons with $n=11$.
To specify a twisted 11-gon, we need 11 vectors in $\Q^3$
(which we then project to 11 vertices in $\Q\PP^2$) and a monodromy from
$SL_3$, which can be defined by fixing 3 more vectors  in $\Q^3$.
Overall we choose 14  vectors in $\Q^3$  uniformly distributed in $[1,10]^3$).

\begin{figure}[h!]
\begin{center}
\includegraphics[width=6.5cm]{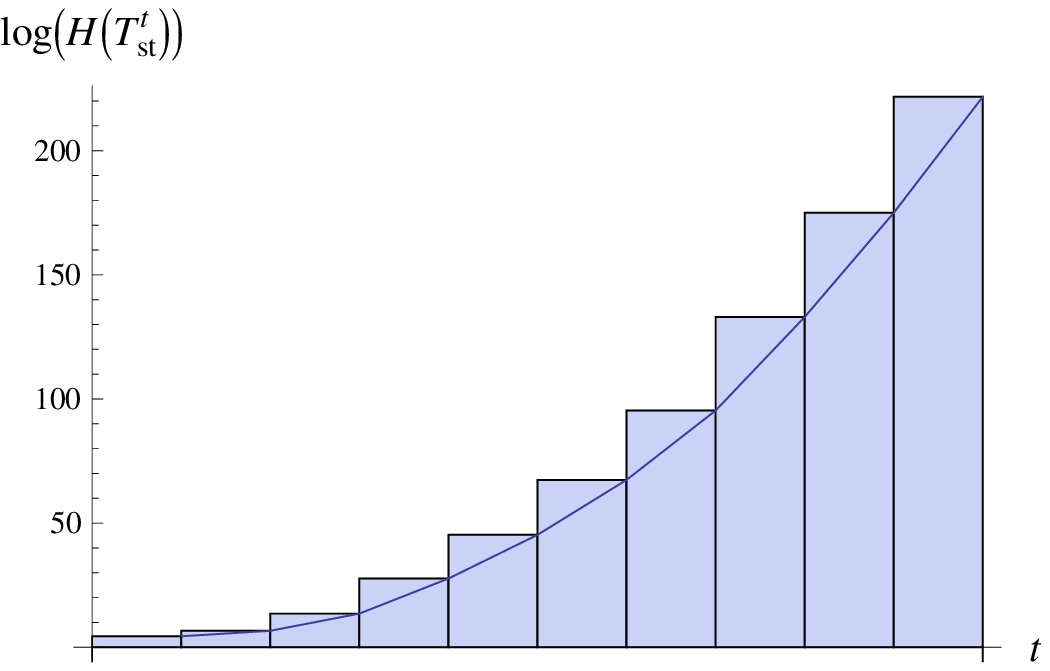}
\hspace{\fill}
\includegraphics[width=6.5cm]{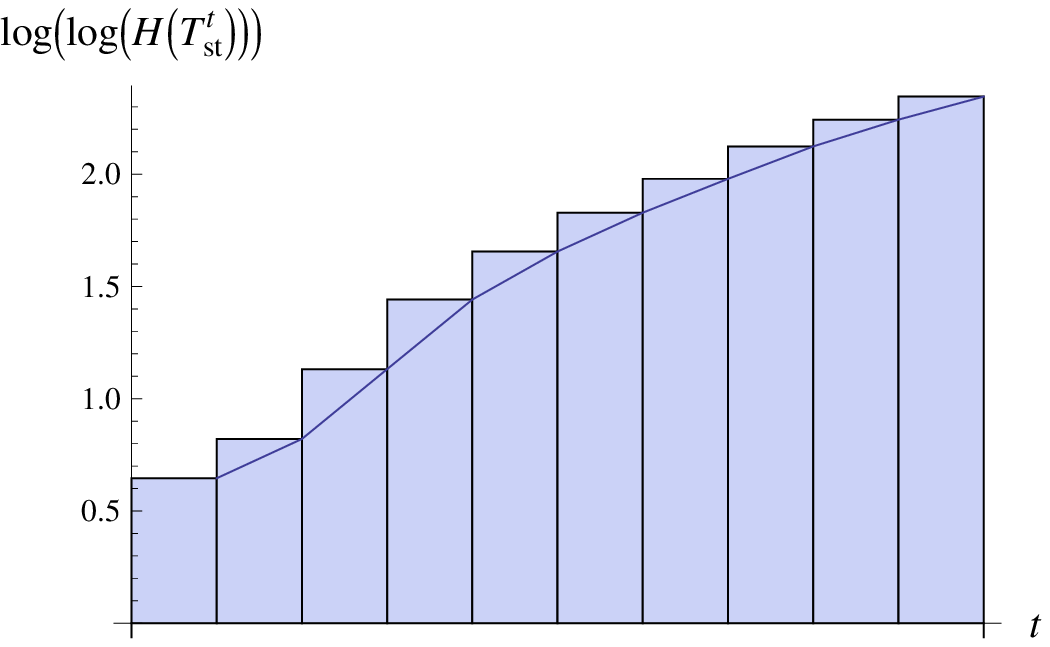}
\caption{Polynomial growth of $\log{H}$ for the map $T_{\rm st}$ in 2D as a function of $t$.}\label{fig:Tst_log}
\end{center}
\end{figure}

First, we start with the standard 2D pentagram map $T_{\rm st}$.
After 10 iterations, the height becomes of the order of $10^{250}$.
Because of its  magnitude, it is natural to use the $\log$ scale and
even the $\log$-$\log$ scale for the height, see Figure~\ref{fig:Tst_log}.

More generally, we are going to study the following $T_{I,J}$ maps in
2D with $I=(i)$ and $J=(j)$, where the
diagonals are chords $P_k=(v_k, v_{k+i})$ and the maps are defined by
intersecting those chords: $T_{I,J}v_k:=P_k\cap P_{k+j}$.

\begin{figure}[h!]
\begin{center}
\includegraphics[width=6in]{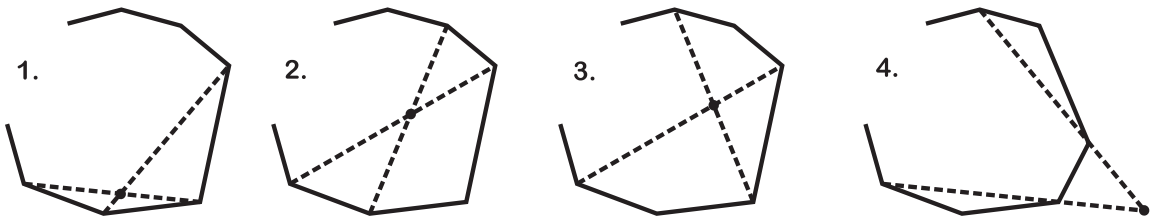}

\caption{The pentagram maps $T_{(i),(j)}$ on twisted
polygons in 2D, where  $T_{(i),(j)} v_k:= (v_k,v_{k+i})\cap(v_{k+j},v_{k+i+j})$, cf. the table below:
1. $T_{(2),(1)}$, 2. $T_{(3),(1)}$, 3.  $T_{(3),(2)}$, 4. $T_{(2),(3)}$.}
\label{fig:T2D-1}
\end{center}
\end{figure}

Note that all these maps are integrable: the integrability of $T_{\rm{st}}:=T_{(2),(1)}$
and of $T_{(3),(1)}$ follows from \cite{OST99, GSTV} (see also Theorem \ref{th-lax}
in Section \ref{intro}). The integrability in the case of $T_{(3),(2)}$, as
well as of its dual $T_{(2),(3)}=T_{(3),(2)}^{-1}\circ Sh$, on 2D  $n$-gons can be proved in a similar way by
changing numeration of vertices, at least for $n$ mutually prime with
$i$ or $j$ and closed polygons. The integrability of such pentagram maps was also observed
experimentally in the applet of R.~Schwartz (personal communication).

In the table below we collected the order of magnitude for the height
growth after 10 iterations for the following maps, see Figure \ref{fig:T2D-1}:

\begin{center}
\begin{tabular}{||c|c|c|c||}
\hline
\multirow{2}*{\#}    & notation for &\multirow{2}*{definition of $Tv_k$}   & height $H$ after       \\
&   pent. map $T=T_{I,J}$   & &10 iterations \\ \hline
1.&$T_{\rm{st}}=T_{(2),(1)}$ &$(v_k,v_{k+2})\cap(v_{k+1},v_{k+3})$ & $10^{320}$ \\ \hline
2.&$T_{(3),(1)}$&$(v_k,v_{k+3})\cap(v_{k+1},v_{k+4})$ & $10^{350}$    \\ \hline
3.&$T_{(3),(2)}$&$(v_k,v_{k+3})\cap(v_{k+2},v_{k+5})$ & $10^{750}$  \\ \hline
4.&$T_{(2),(3)}$&$(v_k,v_{k+2})\cap(v_{k+3},v_{k+5})$ & $10^{800}$  \\ \hline
\end{tabular}

\end{center}

\bigskip

{\it The skew case.}
In all the cases above the pentagram maps were defined by taking intersections of the same type diagonals
at each step. Now we generalize the definition of the classical 2D pentagram map to allow
intersection of {\it different type} diagonals at each step.

\begin{figure}[h!]
\begin{center}
\includegraphics[width=7.5cm]{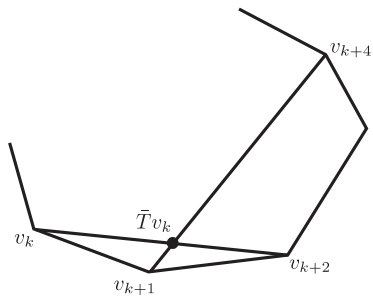}
\caption{The skew pentagram map $\bar T$ is obtained by intersecting diagonals
of length 2 and 3 at each step.}
\label{fig:Tbar_2D}
\end{center}
\end{figure}

As an example, we define the {\it skew pentagram map}  on twisted
polygons in $\PP^2$  by intersecting at each step a short diagonal
$(v_k, v_{k+2})$ of ``length'' 2 and a longer diagonal
$(v_{k+1}, v_{k+4})$ of ``length'' 3:
$\bar Tv_k:=(v_k, v_{k+2})\cap (v_{k+1}, v_{k+4})$, see Figure \ref{fig:Tbar_2D}.
(This map can be described as a universal map $T_{\mathcal I,J}$ of Section \ref{sect:univ}, also cf. \cite{Beffa}, where
for any $d$ one intersects a short diagonal with a hyperplane.)
Note that the skew map $\bar T$ is not a  generalized map of type $T_{I,J}$
from Definition \ref{def:IJ} for any tuples $I$ and $J$, as the latter maps were using the same definition of diagonals at each step.
Now we are going to compare the height growth  for this map
$\bar T$, as well as for several
similar maps, with that for the previously discussed integrable pentagram maps in 2D.

\begin{figure}[h!]
\begin{center}
\includegraphics[width=6.5cm]{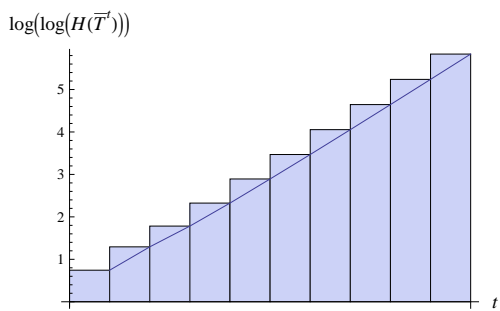}
\caption{Linear growth of $\log{\log{H}}$ for the skew pentagram map $\bar T$ in 2D, which indicates super fast growth
of its height.}\label{fig:Tni_loglog}
\end{center}
\end{figure}

\smallskip

\begin{figure}[h!]
\begin{center}
\includegraphics[width=6in]{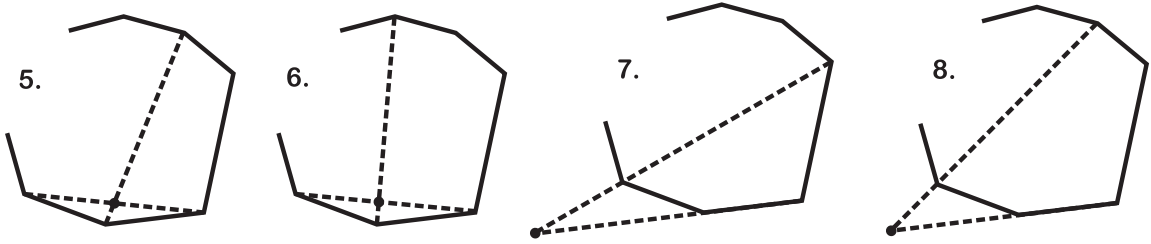}

\caption{The 2D pentagram maps with the following diagonals  intersecting:
5. $(v_k,v_{k+2})\cap(v_{k+1},v_{k+4})$,
6. $(v_k,v_{k+2})\cap(v_{k+1},v_{k+5})$,
7.  $(v_{k+1},v_{k+2})\cap(v_k,v_{k+3})$,
8. $(v_{k+1},v_{k+2})\cap(v_k,v_{k+4})$,
cf. the table above.}\label{fig:T2D-2}
\end{center}
\end{figure}

It turns out that for the skew map $\bar T$ after 10th iteration the
height reaches the order of $10^{10^6}$, see
Figure~\ref{fig:Tni_loglog}. The same order of magnitude for the
height growth is observed for several similar maps,
as summarized in the following table below. We sketch the corresponding
diagonals for these maps on Figure \ref{fig:T2D-2}.

\begin{center}
\begin{tabular}{||c|c|c||}
\hline
\multirow{2}*{\#}   &\multirow{2}*{definition of $Tv_k$}     & height $H$ after       \\
&   &10 iterations \\ \hline
5.&$\bar Tv_k:=(v_k,v_{k+2})\cap(v_{k+1},v_{k+4})$ & $10^{10^6}$  \\ \hline
6.& $(v_k,v_{k+2})\cap(v_{k+1},v_{k+5})$ & $10^{10^6}$    \\ \hline
7.& $(v_{k+1},v_{k+2})\cap(v_k,v_{k+3})$ & $10^{10^6}$   \\ \hline
8.&$(v_{k+1},v_{k+2})\cap(v_k,v_{k+4})$ & $10^{10^6}$  \\ \hline
\end{tabular}

\end{center}

Such a super fast growth is in sharp contrast with the classical
integrable cases
discussed earlier and suggests nonintegrability of all these skew pentagram maps.

\begin{remark}
Note that the above classical and skew examples are conjectured to be
Hamiltonian, regardless of whether they are integrable or not, see
Conjecture \ref{conj:uni} $a)$ and \cite{KS13}.  An example of a different type, the projective heat map in 2D,  was proposed in \cite{heat}: it can be thought of as a dissipative system
on polygons, while its
continuous analog is the curvature flow on curves. This map turns out
to converge to a (projectively) regular $n$-gon, at least for $n=5$.
Such a dynamical system cannot be integrable due to ``dissipation,"
and this non-integrability is of ``non-Hamiltonian" nature.
\end{remark}

\bigskip


\section{Non-integrability in 3D}\label{nonint3D}

In this section we present the results of the numerical integrability test
for various 3D pentagram maps.
First of all,  note that the definition of the height can be
naturally extended to twisted rational polygons in any dimension. For instance, in 3D we employ the cross-ratio coordinates $(x,y,z)$ (defined in \cite{KS}) on the space of twisted $n$-gons ${\mathcal P}_n$ in $\PP^3$, with rational coordinates.\footnote{In any dimension one may use the quasi-periodic coordinates to construct cross-ratio-type coordinates, see  \cite{KS13}.}

\begin{definition}
{\rm
The {\it height} of a twisted $n$-gon $P \in {\mathcal P}_n$ in $\Q\PP^3$ is
$$
H(P):=\underset{0 \le i \le n-1}{\text{max}} \; \text{max}(ht(x_i),ht(y_i), ht(z_i)).
$$
}
\end{definition}

Similarly to the above analysis we trace how fast the height of an initial $n$-gon   for $n=11$ grows
with the number of iterates of different pentagram maps in 3D.
Now we specify 15=11+4  vectors in $\Q^4$ to fix a twisted  $11$-gon in $\Q\PP^3$
and its monodromy from $SL_4$. Again, their coordinates are randomly distributed in $[1,10]$.

It turns out that in 3D there also exists a sharp contrast in the
height growth for different maps. However, the borderline between integrable and non-integrable ones does not lie between the classical and skew cases, and it is more difficult to describe. This is why we group the numerically integrable and non-integrable cases separately.

\bigskip

\begin{figure}[h!]
\begin{center}
\includegraphics[width=6.5cm]{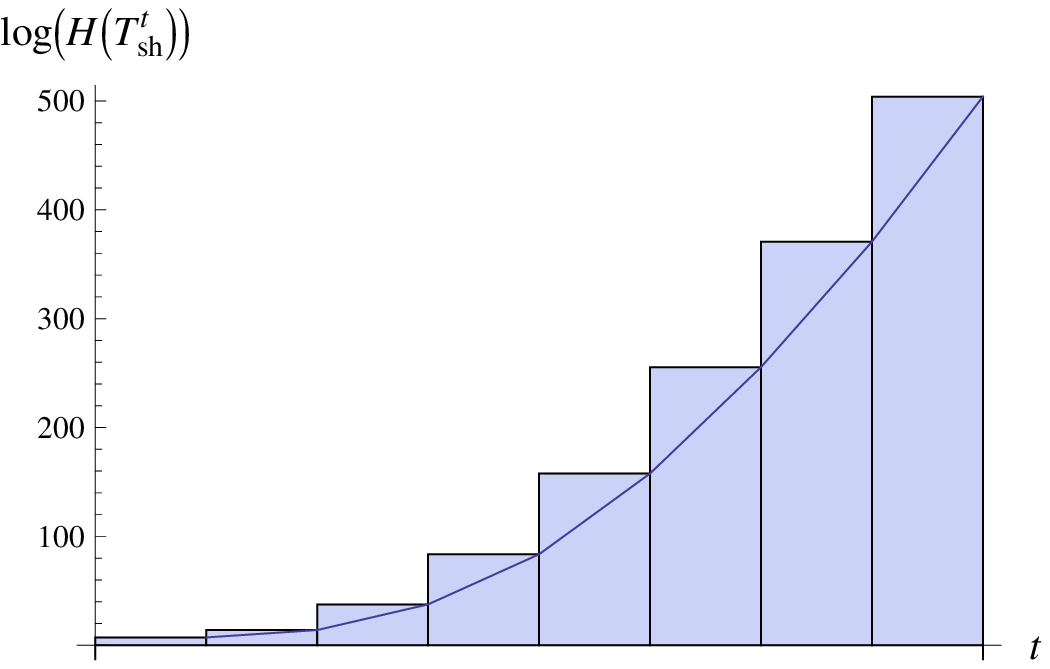}
\hspace{\fill}
\includegraphics[width=6.5cm]{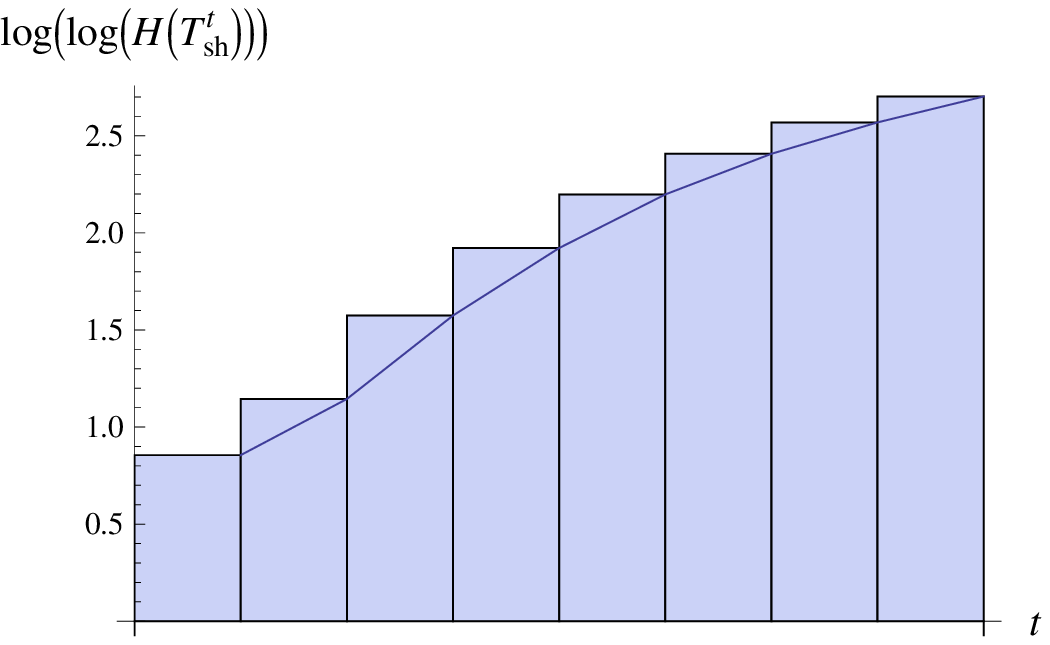}
\caption{Polynomial growth of $\log{H}$ for the integrable pentagram map $T_{\rm sh}$ in 3D
as a function of $t$.}\label{fig:T21_log}
\end{center}
\end{figure}

{\it Numerically integrable 3D cases.}
We start this study with the short-diagonal map $T_{\rm sh}$ in 3D, which is known to be (algebraic-geometric) integrable \cite{KS}. After 8 iterations  of this map, the height of the twisted
$11$-gon in $\Q\PP^3$ becomes of the order of $10^{500}$, see Figure \ref{fig:T21_log}.

The height also grows moderately fast for another integrable map, dented map $T_1$,
reaching the value of the order of $10^{800}$. We also observe a similar moderate growth
for the  (integrable) deep-dented map $T_{m,p}$ in 3D with $m=1$ and $p=3$:
the height remains around $10^{1000}$.

The above cases correspond to taking intersections of consecutive planes, i.e., to $T_{I,J}$ with $J=(1,1)$.
One can also observe a moderate height growth for several pentagram maps with more elaborate tuples $J$.
The results are collected in the following table:

\begin{center}
\begin{tabular}{||c|c|c|c||}
\hline
\multirow{2}*{\#}   & notation for &\multirow{2}*{$Tv_k$ in 3D case}   & height $H$ after\\
  &               pent. map $T=T_{I,J}$  &                           & 8 iterations \\ \hline
1.& $T_{\rm sh}:=T_{(2,2), (1,1)}$  &
$P_k=(v_k, v_{k+2}, v_{k+4})$,  $ Tv_k:= P_k\cap  P_{k+1}\cap  P_{k+2} $  & $10^{500}$\\ \hline
2.& $T_1 :=T_{(2,1), (1,1)}$   &
$P_k=(v_k, v_{k+2}, v_{k+3})$,  $ Tv_k:= P_k\cap  P_{k+1}\cap  P_{k+2} $  & $10^{800}$\\ \hline
3.& $T_{1,3}:=T_{(3,1), (1,1)}$  &
$P_k=(v_k, v_{k+3}, v_{k+4})$,  $ Tv_k:= P_k\cap  P_{k+1}\cap  P_{k+2} $   & $10^{1000}$\\ \hline
4.& $T_{(2,2), (1,2)}$   &
$P_k=(v_k, v_{k+2}, v_{k+4})$,  $ Tv_k:= P_k\cap  P_{k+1}\cap  P_{k+3} $   & $10^{1000}$\\ \hline
5.& $T_{(1,2), (1,2)}$   &
$P_k=(v_k, v_{k+1}, v_{k+3})$,  $ Tv_k:= P_k\cap  P_{k+1}\cap  P_{k+3} $   & $10^{2000}$\\ \hline
6.& $T_{(1,3), (1,3)}$   &
$P_k=(v_k, v_{k+1}, v_{k+4})$,  $ Tv_k:= P_k\cap  P_{k+1}\cap  P_{k+4} $     & $10^{3000}$\\ \hline
\end{tabular}

\end{center}

\begin{remark}\label{int-sum}
The first three cases in the table, with $J=(1,1)$, have been proved to be integrable.
The integrability of the other three, with non-unit $J$, is unknown.
Also, the pattern, which differs these cases from the non-integrable
ones discussed below  is yet to be established.

Note that the case  $T_{(2,2), (1,2)}$ has the ``classical"
diagonal $P_k=(v_k, v_{k+2}, v_{k+4})$, and the pentagram map would be integrable for $J=(1,1)$. However, here for $J=(1,2)$ we  take the intersection
of two consecutive diagonals $P_k$'s and one apart.
(Similarly behaves the map $T_{(2,1), (2,2)}$, which is inverse of $T_{(2,2), (1,2)}$, and hence has the same integrability properties.)
For the last two cases 5) and 6) in the table, their diagonals with $ I=(1,2)$ and $I=(1,3)$
are also  known to be integrable in combination with $J=(1,1)$, while now the numerical results
show that they are also integrable in combination with $J=(1,2)$ and $J=(1,3)$, respectively.
(Note that $T_{I,I^*}$ is always the identity map modulo a shift of
indices,  as follows from the properties of the duality maps, see \cite{KS13}. In particular, e.g., one has $T_{(1,2), (2,1)}=Sh$.)

However, as we will see below, the same diagonals in combination with other $J$'s
may give non-integrability: see the cases 7) and 8) in the table below, where
one mixes $I=(1,2)$ with $J=(3,1)$ or $J=(1,3)$. (Due to duality  \eqref{eq:duality}  one can interchange
these $I$ and $J$, which would lead to the same result on numerical
non-integrability for  $I=(1,3)$ and $J=(2,1)$.)
\end{remark}

\bigskip


{\it Numerically non-integrable 3D cases.}
Non-integrability of pentagram maps appears in several different situations.
As we mentioned in Remark \ref{int-sum}, it can be obtained by taking an ``unusual" intersection tuple
$J$ with a ``usually integrable" jump tuple $I$.

\begin{figure}[h!]
\begin{center}
\includegraphics[width=4.5cm]{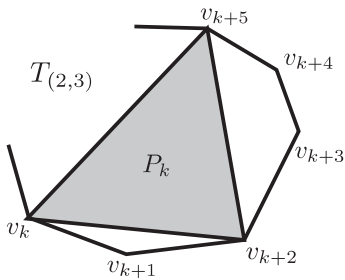}
\caption{The pentagram map $T_{(2,3)}$ in 3D is defined by intersecting three diagonals
$P_k^{(2,3)}:=(v_k, v_{k+2}, v_{k+5})$.}\label{fig:T23_3D}
\end{center}
\end{figure}

Another way to observe non-integrability is to choose a jump tuple $I$
not covered by the integrability theorems (see the survey in Section \ref{intro}).
In 3D we proved integrability for pentagram maps defined by hyperplanes $P_k$
of the following types: $P_k^{\rm st}=(v_k, v_{k+2}, v_{k+4})$, $P_k^1=(v_k, v_{k+2}, v_{k+3})$,
$P_k^{1,p}=(v_k, v_{k+p}, v_{k+p+1})$ (and similarly for $P_k^2$ and $P_k^{2,p}$).
One of the first cases not covered by these results is the pentagram map $T_{(2,3)}:=T_{(2,3),(1,1)}$ in 3D
defined by the hyperplanes
$P_k^{(2,3)}:=(v_k, v_{k+2}, v_{k+5})$ with the jump tuple $I=(2,3)$, while $J=(1,1)$, in notations of \cite{KS13}:
$$
T_{(2,3)}v_k:= (v_k, v_{k+2}, v_{k+5}) \cap (v_{k+1}, v_{k+3}, v_{k+6})
 \cap (v_{k+2}, v_{k+4}, v_{k+7}),
$$
see Figure \ref{fig:T23_3D}.
We conjectured in \cite{KS13} (see also Conjecture \ref{conj:uni} $a)$  that all maps defined by taking intersections of the same
diagonals are  discrete Hamiltonian. But they still might be non-integrable and $T_{(2,3)}$
is the first candidate for that. Here we present a numerical evidence for such a non-integrability.

\begin{figure}[h!]
\begin{center}
\includegraphics[width=6.5cm]{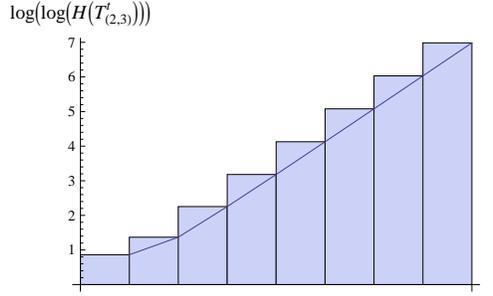}
\caption{Linear growth of $\log{\log{H}}$ for the map $T_{(2,3)}$  in 3D indicating super fast growth
of its height and apparent non-integrability.}
\label{fig-nonint}
\end{center}
\end{figure}

The height growth turns out to be enormously faster for the map $T_{(2,3)}$ than for all integrable maps discussed above:
after 8 iterations the height already reaches the order of magnitude
of over $10^{10^7}$, see Figure \ref{fig-nonint}.
The map $T_{(2,3)}$ in 3D is a map defined by the same diagonal plane at each step, i.e., it is
of type $T_{I,J}$. However,  in a sense it is mimicking the skew map $\bar T$ in 2D defined by different type diagonals.
More cases of presumably non-integrable maps are given in the table below.
We  mention that the case 11) was discussed in \cite{KS},
where the problem of its integrability was posed. It looked
conjecturally integrable as the corresponding pentagram map $T_{(3,3), (1,1)}$ is the intersection of three consecutive very symmetric diagonals $P_k:=(v_k, v_{k+3}, v_{k+6})$. However, the current numerical evidence suggests its non-integrability.

Finally, one more source of would-be non-integrable maps are skew pentagram maps, and, in particular, the maps constructed by intersecting different planes of complimentary dimensions
at each step. The  cases 12) and 13) in the table below illustrate the latter.

\begin{center}
\begin{tabular}{||c|c|c||}
\hline
\multirow{2}*{\#} &  \multirow{2}*{$Tv_k$ in 3D case}            & height $H$ after\\
                                &                               & 8 iterations \\ \hline
7.&$T_{(1,2), (3,1)}$                                 & $10^{3 \cdot 10^7}$\\ \hline
8.&$T_{(1,2), (1,3)}$                                 & $10^{3 \cdot 10^7}$\\ \hline
9.&$T_{(2,3)}:=T_{(2,3), (1,1)}$                                 & $10^{10^7}$\\ \hline
10.&$T_{(2,4)}:=T_{(2,4), (1,1)}$                                 & $10^{10^7}$\\ \hline
11.&$T_{(3,3)}:=T_{(3,3), (1,1)}$                                 & $10^{10^7}$\\ \hline
12.&$(v_k,v_{k+3})\cap(v_{k+1},v_{k+2},v_{k+4})$& $10^{10^6}$\\ \hline
13.&$(v_{k+1},v_{k+3})\cap(v_k,v_{k+2},v_{k+5})$& $10^{10^6}$\\ \hline
\end{tabular}

\end{center}

\begin{remark} Note that all maps considered above can be expressed by rational functions.
Apparently, in  integrable cases many cancellations of different terms
occur after several iterations,  resulting in a much lower complexity,
while one does not have those cancelations in   non-integrable cases.
It would be very interesting to prove the observed non-integrability directly, rather than numerically, possibly by employing  Ziglin's or Morales--Ramis' methods.
\end{remark}

The above consideration provides an evidence for Conjecture~\ref{conj:uni}~$b)$:
A necessary condition for integrability of a  pentagram map is to be of $T_{I,J}$-type,
i.e., to be defined by intersections of the same-type diagonal hyperplanes at each step
(see Definition \ref{def:IJ}).
While the class of universal  pentagram maps $T_{\mathcal I,J}$ is
very broad,  all known {\it integrable} examples (such as short-diagonal, dented,
deep-dented, corrugated and partially corrugated pentagram maps) can  indeed be
presented as  $T_{I,J}$-type maps.
However, this condition is not sufficient for integrability, as many examples of this section indicate.

\bigskip

{\bf  Acknowledgments.}
B.K. is grateful to the Simons Center for Geometry and Physics for support and kind hospitality.
The research of B.K. and F.S. was partially supported by NSERC grants.

\end{document}